\documentclass[a4paper]{article}
\usepackage{graphicx} 
\usepackage{amsmath, amsthm, amsfonts, amssymb, amscd, bm, enumerate}
\usepackage{verbatim}

\usepackage{color}
\usepackage{tikz}

\usepackage{url}

\renewenvironment{proof}{{\bfseries Proof:}}{\qed}
\newtheorem{theorem}{Theorem}
\newtheorem{proposition}{Proposition}

\newtheorem{lemma}{Lemma}
\newtheorem{definition}{Definition}
\newtheorem{assumption}{Assumption}
\newtheorem{remark}{Remark}

\DeclareMathOperator{\Tr}{Tr}

\def\B{\mathcal{B}}
\def\F{\mathcal{F}}

\def\L{\mathcal{L}}

\def\U{\mathcal{U}}

\def\bR{\mathbb{R}}
\def\bE{\mathbb{E}}
\def\bN{\mathbb{N}}

\def\bone{\mathbf{1}}
\def\bP{\mathbb{P}}

\def\bQ{\mathbb{Q}}

\mathchardef\mhyphen="2D

\def\tT{\tilde{T}}

\def\u{\bar{u}}
\def\tu{\tilde{u}}
\def\tx{\tilde{x}}

\def\X{\bar{X}}
\def\Y{\bar{Y}}
\def\Z{\bar{Z}}
\def\<{<}
\def\>{>}

\begin{document}
\title{Classical Adjoints for Ergodic Stochastic Control}
\author{Samuel N. Cohen\\ University of Oxford \\ \\ Victor Fedyashov\\ University of Oxford}

\date{\today}

\maketitle

\begin{abstract}

In this paper we consider ergodic optimal control of a diffusion process $\{X^u_t\}_{t \geq 0}$, taking values in $\bR^n$, where both drift and volatility are controlled. We establish a novel strong duality between the existence of a unique solution to the infinite horizon adjoint BSDE and strong dissipativity of $X^u$. We then proceed to show that the latter implies irreducibility, the strong Feller property and exponential ergodicity. We conclude by discussing the connection with ergodic BSDEs. 



\end{abstract}

\section{Introduction}

Ergodic stochastic optimal control is a growing area of optimal control theory that is trying to understand optimisation with an average cost criterion. In other words, it is concerned with payoffs that value the future as much as the present. So far, the main probabilistic framework that has been introduced to address the problem is ergodic BSDEs (see, e.g. \cite{Levy_paper}, \cite{Hu}, \cite{Hu_Banach}). However, this can only be applied to the so called weak control formulation: every control is treated as generating a change of measure, so the value functional takes the form 
\[
	J(x, u) = {\lim\sup}_{T\to\infty} T^{-1} \bE^{u}\bigg[\int_0^T L(X_t, u_t) dt\bigg],
\]
where $X$ represents the forward process, and the control $\{u_t\}_{t \geq 0}$ is an $\F_t$-predictable process taking values in a separable locally compact metric space which determines the expectation $\bE^u$. It is well known that, if we allow for control of both drift and volatility, such representation becomes impossible. In order to deal with these problems in their strong formulation, numerous stochastic maximum principles have been established in both finite (first introduced in \cite{Peng_SMP}) and recently infinite (see \cite{Oks}) horizons. In this paper we propose a candidate for the adjoint BSDE in the case of ergodic control and explain why it is a natural one. We therefore conclude that existence of a solution to this BSDE is necessary for optimality. Assuming convexity of the Hamiltonian, we then prove a version of the sufficient stochastic maximum principle. We also observe that there is a strong connection between the adjoint equation admitting a solution and the ergodicity of the forward process. 

The rest of the paper is organised as follows: in the next section we setup the problem. Section 3 deals with stochastic maximum principle approach, and its connections with strong dissipativity of the forward process under all controls. Section 4 shows that strong dissipativity, in our general setting, implies exponential ergodicity, establishing irreducibility and strong Feller properties along the way. In Section 5 we look at an EBSDE representation of the problem. Section 6 concludes.

\section{Setup}

In this section we introduce the preliminary concepts and definitions we will be using repeatedly in the sequel. We also state the principle problem of interest and discuss its connections to the existing theory. 

\subsection{The forward dynamics}

We begin by introducing the forward process $\{X^u_t\}_{t \geq 0}$. Let the controlled dynamics of $\{X^u_t\}_{t \geq 0}$ be governed by the It\^o SDE:
\begin{equation}
	dX^u_t = b(t,X^u_t,u_t)dt + \sigma(t,X^u_t,u_t)dW_t, \quad X_0 = x_0 \in \bR^n,
\label{CDX}
\end{equation}
where $W$ denotes a standard Brownian motion under the measure $\bP$. We need the following assumptions:
\begin{assumption} Conditions to guarantee existence of a strong solution to (\ref{CDX}):
\begin{itemize}
\item The triple $(\Omega, \F, \bP)$ is a complete probability space, $(\F_t)_{t\geq 0}$ is an augmented filtration satisfying the usual conditions, and $(W_t)_{t \geq 0}$ is a Brownian motion. 
\item The coefficients $b(x,t,u)$ and $\sigma(x,t,u)$ are measurable in $u$ and $t$, locally Lipschitz in $x$, that is, for every $T,N$, there exist a constant $K$ depending only on $T$ and $N$, such that
\[
	| b(t,x,u) - b(t,y,u)| + |\sigma(t,x,u)-\sigma(t,y,u)| < K |x-y|
\]
holds for all $u \in \U$, $|x|,|y| \leq N$ and all $t \leq T$. 
\item The linear growth condition holds. In other words, for all $u \in \U, x \in \bR^n, t \geq 0$,
\[
	|b(t,x,u)| + |\sigma(t,x,u)| \leq \bar{K}(1 + |x|),
\]
where $\bar{K}$ is a constant. 
\item $x_0 \in \bR^n$ is a constant. 
\end{itemize}
\label{existence_solution}
\end{assumption}
\noindent The first part of the following result is standard. The second is an easy application of Gr\"onwall's lemma (for details see, e.g. Theorem 6.30 in \cite{Fima}):
\begin{theorem} If Assumption \ref{existence_solution} is satisfied then, for any fixed control process $\{u_t\}_{t \geq 0}$, there exists a unique strong solution to the equation (\ref{CDX}). Moreover,
\[
	\bE \bigg( \sup_{0 \leq t \leq T} | X^u_t |^2 \bigg) < C (1 + \bE|x_0|^2),
\]
where the constant $C$ depends only on $K$ and $T$.
\end{theorem}
\noindent For the rest of the paper we also make the following assumptions: 
\begin{assumption}\begin{enumerate}[(i)]
\item There exist constants $\underline{\sigma}$ and $\overline{\sigma}$, such that
\[
	\underline{\sigma} \leq \| \sigma(t,x,u) \| + \| \sigma^{-1}(t,x,u) \| \leq \overline{\sigma}
\]
holds for all $u \in \U,x \in \bR^n,t \geq 0$.
\item There exists a constant $\tilde{C} > 0$, such that
\[
	\| b(t,0,u) \| \leq \tilde{C}
\]
holds for all $u \in \U, t \geq 0$.
\end{enumerate}
\label{drift}
\end{assumption}

\begin{remark} Note that Lipschitz continuity of the coefficients is stronger than what is generally required for the purposes of control. However, in the sequel we will need certain estimates for the stability of the forward process in its initial value, and stronger regularity assumptions will come into play. 
\end{remark}

\noindent The following definition will prove crucial in the sequel, since it is very closely related to ergodicity of the solution process to (\ref{CDX}).
\begin{definition} Let the dynamics of $\{X^u_t\}_{t \geq 0}$ be governed by (\ref{CDX}). We say that $X^u$ has a \textit{strongly dissipative} drift (or $X^u$ is strongly dissipative), if there exists a constant $\mu > 0$ such that
\[
	\langle  b(t, x,u_t) - b(t,y,u_t), x-y   \rangle \leq - \mu \| x - y \|^2
\]
holds for all $t \geq 0$. 
\label{SDD}
\end{definition}

\subsection{The problem}

\begin{definition} Let $\U$ be a separable metric space. Then a stochastic process $\{u_t\}_{t \geq 0}$ is called admissible if it is predictable in the filtration $\{\F_t\}_{t \geq 0}$ with values in $\U$.
\end{definition}

\noindent We consider the ergodic control problem, namely that of minimising (over the space of admissible controls) the functional 
\begin{equation}
	J(x_0, u) = {\lim\sup}_{T\to\infty} T^{-1} \bE \bigg[\int_0^T L(t,X^u_t, u_t) dt\bigg],
\label{F1}
\end{equation}
where the dynamics of the controlled process $\{ X^u_t \}_{t \geq 0}$ are given by (\ref{CDX}), and the coefficients $b(\cdot,u)$ and $\sigma(\cdot,u)$ satisfy Assumptions \ref{existence_solution} and \ref{drift} for all admissible controls $\{u_t\}_{t \geq 0}$. The finite horizon case, namely the problem of minimising 
\begin{equation}
	J(x_0, u) =  \bE \bigg[\int_0^T L(t,X^u_t, u_t) dt + g(X^u_T)\bigg], 	
\label{FHJ}
\end{equation}
is well understood (see, e.g. the seminal paper by Peng \cite{Peng_SMP}). The adjoint equation is 
\begin{equation}
	dY(t) = - \nabla_x H(t,X^u_t,u_t,Y(t),Z(t))dt + Z(t)dW_t, \quad Y(T) = \nabla_x g(X_T),
\label{adjoint_BSDE}
\end{equation}
where $H:[0,T]\times \bR^n \times \U \times \bR^n \times \bR^{n \times n} \to \bR$ is the Hamiltonian, defined by 
\begin{equation}
	H(t,x,u,y,z) := b'(t,x,u)y + \Tr(\sigma'(t,x,u)z) + L(t,x,u).
\label{Ham}
\end{equation}
and it is assumed that $H$ is differentiable in $x$ (\textit{NB: here and for the rest of the paper $'$ denotes the transpose}). The following condition for the necessary stochastic maximum principle can then be established (see Theorem 3 in \cite{Peng_SMP}):
\begin{theorem} Suppose $\{u_t\}_{t \geq 0}$ is the optimal control process, and $\{X^{u}_t\}_{t \geq 0}$ are the corresponding controlled dynamics for the problem (\ref{FHJ}). Then there exists a solution to the BSDE (\ref{adjoint_BSDE}) on finite horizons.
\label{pengSMP}
\end{theorem}

\begin{remark} We note that, in \cite{Peng_SMP}, since no convexity on the Hamiltonian is assumed, an additional second order adjoint equation appears. However in this paper we will focus on the first order BSDE. 
\end{remark}

\noindent One can also show the connection between the adjoint process $\{Y(t) \}_{t \geq 0}$ and the dynamic programming principle. The following result can be found, for example, in Chapter 6 of \cite{Pham}:

\begin{lemma} Let $\{X^{u,t,x}_s\}$ denote the solution to (\ref{CDX}) with the initial condition $X^{u,t,x}_t = x$. Define the value function $v(\cdot,\cdot)$ as follows:
\[
	v(t,x) := \inf_{u} \bE \bigg[ \int_t^T L(s,X^{u,t,x}_s,u_s)ds + g(X^{u,t,x}_T) \bigg].
\]
Then 
\[
	\bar{Y}(t) = \nabla_x v(t,\bar{X}_t),
\]
where $\bar{Y}$ is the first component of the solution to the adjoint BSDE with optimal control process $\{ \bar{u}_t \}_{t\in [0,T]}$ and the corresponding controlled forward process $\{\bar{X}\}_{t \geq 0}$.
\label{DPP}
\end{lemma}
%%%%%%%
\begin{remark}From Lemma \ref{DPP}, we see that the process $\bar{Y}$ can be thought of as the optimal marginal value of a change in $X$. This will be an intuitive starting point when developing an infinite horizon analogue for the adjoint BSDE in the next section.
\label{rem1}
\end{remark}
%%%%%%%

\section{Stochastic maximum principle approach to the problem of ergodic optimal control}

The aim of this section is to extend the finite horizon framework to the ergodic case, showing an analogue of the adjoint equation, and establishing existence and uniqueness of its solutions. For the rest of the paper, we make the following assumption:
\begin{assumption} The space of controls $\U$ is a locally convex real topological vector space. (For the sake of simplicity we consider $\U = \bR^k$, but the subsequent analysis does not change provided that the G\^ateaux derivative is well defined.)
\label{assU}
\end{assumption}

%%%%%%%%%%%%%%

\subsection{Adjoint BSDE in infinite horizon}

In this subsection we explain (\textit{in a heuristic way}) why the adjoint equation to the ergodic control problem becomes an infinite horizon BSDE. For the remainder of this subsection, we denote the optimal control pair as $(\u,\X)$. We first state (without a proof) a proposition that can be easily verified:

\begin{proposition} The following statements about infinite horizon optimisation are equivalent:
\begin{itemize} \item The control process $\{\u_t\}_{t \geq 0}$ minimises the functional $J(x_0,\cdot)$ as defined in (\ref{F1}).
\item The control process $\{\u_t\}_{t \geq 0}$ formally minimises the (non averaged) functional defined by 
\begin{equation}
		\bar{J}(x_0,\cdot) = \bE \bigg[ \int_0^{\infty} L(t,X^{\cdot}_t,\cdot) dt\bigg]
\label{fTmp}
\end{equation}
in the sense that (since the right hand side in (\ref{fTmp}) may diverge) we have 
\[
	\limsup_{T \to \infty} \frac{\bE \big[ \int_0^{T} L(t,X^{\u}_t,\u_t) dt\big]}{\bE \big[ \int_0^{T} L(t,X^{u}_t,u_t) dt\big]} \leq 1
\] 
for all admissible controls $u$. 
\end{itemize}
\label{P1}
\end{proposition}

\noindent The problem of minimising (\ref{fTmp}) is in some sense simpler than the ergodic one. This is because, considered formally, optimality in (\ref{fTmp}) implies optimality over every finite horizon by the dynamic programming principle, whereas the ergodic problem only considers optimality in the limit (and so the dynamic programming principle is not necessary for optimality, see Theorem 5.1 in \cite{Andrew} for an example of this phenomenon). Similarly to Lemma \ref{DPP}, suppose we can define 
\[
	v(t,x) := \inf_{u} \bE \bigg[ \int_t^{\infty} L(s,X^{u,t,x}_s,u_s)ds + g(X^{u,t,x}_T) \bigg].	
\]
We can directly employ the dynamic programming principle to see that, for any $T \geq 0$, we have a finite horizon problem of finding 
\begin{equation}
	\inf_{u}\bar{J}(x_0,u) =  \inf_{u}\bE \bigg[ \int_0^T L(t,X^u_t,u_t) dt\bigg] + v(T,X^u_T).
\label{UU}
\end{equation}
Then, since $\{\u\}_{t \geq 0}$ minimises $\bar{J}(x_0,\cdot)$, it is also optimal for (\ref{UU}) for every $T > 0$. Therefore, by Theorem \ref{pengSMP}, we infer that the equation 
\begin{equation}
	Y(t) = Y(T) + \int_t^T \nabla_x H(s,X^u_s,u_s,Y(s),Z(s))ds - \int_t^T Z(s)dW_s
\label{ihAdj}
\end{equation}
should hold for all $0 \leq t \leq T < \infty$. This suggests strongly that (\ref{ihAdj}) is a good candidate for the adjoint equation in an ergodic maximum principle. 

Since, unlike in the finite horizon case, there is no terminal condition, we need to say something about the limiting (as $T \to \infty$) behaviour of the `optimal' solution to (\ref{ihAdj}) (solution with the optimal control pair $(\u,\X)$). As pointed out in Remark \ref{rem1}, the process $\bar{Y}$ can be thought of as the marginal value of a change in $X$. Therefore the condition $\Y(T) \to 0$ as $T \to \infty$ (which would be the direct analogue of the finite horizon case) does not seem to be the right one. In the next section we will discuss this in more detail.

%%%%%%%%%%%%%%%%

\subsection{Stochastic Maximum Principle}
In the previous section we have explained why (\ref{ihAdj}) is a reasonable candidate for the adjoint equation associated with the problem of ergodic control. One can therefore expect that the existence of a solution to (\ref{ihAdj}) is necessary for optimality, given dynamic programming should hold. The goal of this section is to show that, under certain additional assumptions, it is also sufficient. The main result is the following theorem (similar in spirit to the result in \cite{Oks}), which is a version of the stochastic maximum principle:

\begin{theorem} Suppose the control process $\{\u_t\}_{t \geq 0}$ is admissible, and denote the associated forward process $\{\X_t\}_{t \geq 0}$. Suppose that the corresponding infinite horizon adjoint BSDE (\ref{ihAdj}) has solution $(\Y,\Z)$. Suppose further that, for any admissible control process $\{u_t\}_{t \geq 0}$, 
\begin{equation}
	\limsup_{T \to \infty} \frac{1}{T}\bE \bigg[ \Y_T (X^u_T - \X_T) \bigg] = 0.
\label{adjCond}
\end{equation}
Moreover, assume that the Hamiltonian $H$ as defined in (\ref{Ham}) is convex in $x,u$, differentiable with respect to $x,u$, and 
\begin{equation}
	  H(t,\X_t,\u_t,\Y(t),\Z(t))  = \min_{u \in \U}  H(t,\X_t,u,\Y(t),\Z(t)).
\label{maxHu}
\end{equation}
Then the control $\{\u_t\}_{t \geq 0}$ is optimal. 
\end{theorem}

\begin{remark} All conditions of this theorem are intuitive apart from (\ref{adjCond}). However, it also has an economic interpretation. It tells us that if the the optimal marginal value of a change in $X$ is positive for some large time $T$, then $X^u_T - \X_T$, the difference in the state between an arbitrary controlled process and the optimal one, is on average sublinear in time. 
\end{remark}

\noindent \begin{proof} In order to prove optimality of $\u$, we need to show that, for all admissible controls $\{u_t\}_{t \geq 0}$, 
\[
	J(x_0,\u) - J(x_0,u) \leq 0.
\]
Since $\limsup$ is sub additive, we have, for any functions $g,f:\bR^+ \to \bR$,
\[
	\limsup_{T \to \infty} g(T) - \limsup_{T \to \infty} f(T) \leq \limsup_{T \to \infty} (g-f)(T),
\] 
and therefore it is enough to show that 
\begin{equation}
	\limsup_{T\to\infty} \frac{1}{T} \bE \bigg[\int_0^T \big( L(t,\X_t, \u_t) - L(t,X^u_t, u_t) \big) dt \bigg] \leq 0.
\label{limS}
\end{equation}
By definition of the Hamiltonian, we can write 
\[
\begin{split}
	L(t,\X_t, \u_t) - L(t,X^u_t, u_t) &= H(t,\X_t,\u_t,\Y(t),\Z(t)) - H(t,X^u_t,u_t,\Y(t),\Z(t)) \\
	&\quad + \Tr \bigg(( \sigma'(t,X^u_t,u_t) - \sigma'(t,\X_t,\u_t))\Z(t)\bigg) \\
	&\quad +  \bigg( b'(t,X^u_t,u_t) -  b'(t,\X_t,\u_t)\bigg)\Y(t)
\end{split}
\]
By convexity of $H(t,\cdot,\cdot,y,z)$ for all $(t,y,z) \in \bR^+ \times \bR \times \bR^n$, we have 
\begin{equation}
\begin{split}
	H(t,\tx,\tu,y,z) &-H(t,x,u,y,z) \\
	&\leq \nabla_x H(t,\tx,\tu,y,z)' (\tx - x) + \nabla_u H(t,\tx,\tu,y,z)' (\tu - u)
\end{split}
\label{eHam}
\end{equation}
Given (\ref{maxHu}), we know that 
\[
	 \nabla_u H(t,\tx,\tu,y,z)' (\tu - u) \leq 0,
\]
and thus (\ref{eHam}) becomes 
\[
	H(t,\tx,\tu,y,z) -H(t,x,u,y,z) \leq -\nabla_x H(t,\tx,\tu,y,z)' (x - \tx).
\]
We recall that the dynamics of $\{\Y_t\}_{t \geq 0}$ are given by
\[
	d\Y_t = - \nabla_x H(t,\X_t,\u,\Y_t,\Z_t)dt + \Z_t dW_t.
\]
Therefore, applying It\^o's formula to $\Y_T (X^u_T - \X_T)$ (note that $\X_0 = X^u_0$) and substituting the result into (\ref{adjCond}), we arrive at
\[
	\limsup_{T\to\infty} \frac{1}{T} \bE \bigg[\int_0^T \big( L(t,\X_t, \u_t) - L(t,X^u_t, u_t) \big) dt \bigg] \leq \limsup_{T \to \infty} \frac{1}{T}\bE \bigg[ \Y_T (X_T - \X_T) \bigg],
\]
concluding the proof. 

\end{proof}

\subsection{Illustration in one dimension}

\noindent So far we have seen the close relationship between optimality of a control $u$ and the existence of solutions to the adjoint BSDE (\ref{ihAdj}). We now wish to ask, for a given control, can we see that this BSDE admits a solution? Given we are working over infinite horizons, this is a somewhat delicate question. 

In this subsection we consider the simplest case, where the process $\{X^u_t\}_{t \geq 0}$ is one-dimensional, in order to motivate the subsequent discussion. In order to separate $Y$ in the driver, we write 
\[
	f(s,x,u,z) := H_x(s,x,u,y,z) - b_x(s,x,u)y.
\]
The first thing we notice is that the adjoint BSDE 
\begin{equation}
	Y(t) = Y(T) + \int_t^T [f(s,X^u_s,u_s,Z(s)) - b_x(s,X^u_s,u_s)Y(s)]ds - \int_t^T Z(s)dW_s 
\label{1D}
\end{equation}
looks very similar to a `discounted' one, where the role of the discount factor is played by the coefficient $b_x(s,X^u_s,u_s)$. It is known in the theory of infinite horizon discounted BSDEs that, in order to guarantee the existence of a solution to (\ref{1D}), it is reasonable to assume (see, e.g. \cite{Royer}) that there exists a constant $k > 0$, such that 
\[
	b_x(s,x,u) < -k
\]
for all $s \geq 0$, $x \in \bR^n$ and $u \in \U$. The reason is the following: if we allow the `discount' to be positive, then there is no hope for a solution in general. To illustrate why, consider the following equation:
\[
	dY_t =  - rY_tdt + Z_tdW_t.
\]
We immediately see that 
\[	
	Y_0 = \bE \big[ e^{rT}Y_T \big], \text{  for all } T \geq t.
\]
Therefore, if $Y_T$ does not vanish, there is no bounded solution, since $e^{rT}$ explodes as $T \to \infty$. The case where $b_x$ is nonpositive is still an open question (it is conjectured that the solution exists if the discount is negative on a set of positive ergodic measure). The following lemma is trivial, but we state it since its generalisation will play an important role in the sequel.

\begin{lemma} Let $g:\bR \to \bR$ be of class $C^1$ and let $k>0$ be a constant. Then the following conditions are equivalent:
\begin{enumerate}[(i)] \item $g'(x) \leq -k$ for all $x \in \bR$.
\item $(g(x)-g(y))(x-y) \leq -k (x-y)^2$.
\end{enumerate}
\label{lemmaTemp}
\end{lemma}

\noindent \begin{proof} (i)$\to$(ii): By the mean value theorem, we have, for some $\bar{x} \in (x,y)$,
\[
	(g(x)-g(y))(x-y) = g'(\bar{x})(x-y)^2 \leq -k (x-y)^2.
\]
(ii)$\to$(i): Dividing both sides by $(x-y)^2$, we obtain 
\[
	\frac{g(x)-g(y)}{x-y} \leq -k.
\]
Setting $x = y+h$, and passing to the limit, the result follows. 

\end{proof}

We notice that the condition (ii) is nothing else but the strong dissipativity (in the sense of Definition \ref{SDD}) of the drift of the process $\{X^u\}_{t \geq 0}$ for all controls $u$. In the sequel we will establish (see Theorem \ref{mainThm}) that the latter implies ergodicity. Therefore Lemma \ref{lemmaTemp} shows the equivalence between the condition required for the existence of a solution to the adjoint BSDE in infinite horizon, and the one ensuring ergodicity of the forward process under all possible controls. This conclusion is relatively intuitive: given enough assumptions so that we can solve an optimal ergodic control problem, we should know that the controlled process is ergodic! What is surprising is that the converse is also true -- given the dissipativity of the controlled dynamics, one can hope for classical solutions to the infinite horizon adjoint BSDE. 

 The question becomes: what happens in the multidimensional setting? Is it possible to generalise this connection? The rest of the paper shows that the answer is yes.

\subsection{General case}

\noindent In the general multiple dimensional case, the notion of `discounting' becomes nontrivial. Therefore in order to ensure existence and uniqueness of a bounded solution to the multidimensional infinite horizon adjoint BSDE (\ref{IHBSDE}), we impose the following conditions (inspired by \cite{Royer}) on the gradient of the drift of $X^u$:

\if 0

\begin{assumption} 
\begin{itemize}

\item Functions $b(t,\cdot,u)$, $\sigma(t,\cdot,u)$ and $L(t,\cdot,u)$ are twice continuously differentiable.

\item There exists constants $\bar{C},k > 0$, such that
\begin{equation}
	\langle \nabla_x b(t,x,u)y,y \rangle \leq -k \| y \|^2 + \bar{C}
\end{equation}
holds for all $t \geq 0$, $x,y \in \bR^n$ and all $u \in \U$. 
\item  The driver $\nabla_x H(t,x,u,y,z)$ is uniformly Lipschitz in $z$.
\item  The driver $\| \nabla_x H(t,x,u,0,0) \|$ is uniformly bounded by a constant $C \in \bR$. 
\end{itemize}
\label{assumption_gradient}
\end{assumption}

\fi

\begin{assumption} 
\begin{itemize}

\item Functions $b(t,\cdot,u)$, $\sigma(t,\cdot,u)$ and $L(t,\cdot,u)$ are twice continuously differentiable.

\item There exists a constant $k > 0$, such that
\begin{equation}
	\langle \nabla_x b(t,x,u)y,y \rangle \leq -k \| y \|^2 
\label{E3}
\end{equation}
holds for all $t \geq 0$, $x,y \in \bR^n$ and all $u \in \U$. 
\item  The driver $\nabla_x H(t,x,u,y,z)$ is uniformly Lipschitz in $z$.
\item  The driver $\| \nabla_x H(t,x,u,0,0) \|$ is uniformly bounded by a constant $C \in \bR$. In other words, $\| \nabla_x L(t,x,u) \| \leq C$ for all $(t,x,u) \in \bR^+ \times \bR^n \times \U$. 
\end{itemize}
\label{assumption_gradient}
\end{assumption}

\noindent We are now ready to prove the existence and uniqueness of a bounded solution to the infinite horizon adjoint BSDE 
\begin{equation}
	Y(t) = Y(T) + \int_t^T \nabla_x H(s,X^u_s,u_s,Y(s),Z(s))ds - \int_t^T Z(s)dW_s.
\label{IHBSDE}
\end{equation}

\begin{theorem} Suppose Assumption \ref{assumption_gradient} is satisfied. Then for any $u$ there exists an adapted solution $(Y, Z)$, with $Y$ c\`adl\`ag and $Z\in \L^2(W)$ to the infinite horizon equation (\ref{IHBSDE}) for all $0\leq t\leq T<\infty$, satisfying $\|Y(t) \|\leq C'$ for some $C' \in \bR$, and this solution is unique among bounded adapted solutions.

Furthermore, if $(Y^T, Z^T)$ denotes the (unique) adapted square integrable solution to
\begin{equation}\label{FHBSDE}
	Y^T(t) =  \int_t^T \nabla_x H(t,X^u_s,u_s,Y^T(s),Z^T(s))ds - \int_t^T Z^T(s)dW_s,	
\end{equation}
then $\lim_{T\to\infty} Y^T(t) = Y(t)$ a.s., uniformly on compact sets in $t$.

\end{theorem}

\noindent \textbf{Proof:} We start by proving that if a bounded solution exists, it is unique. Suppose we have two bounded solutions $(Y,Z)$ and $(Y',Z')$ to (\ref{IHBSDE}). We denote $\delta Y := Y - Y'$, $\delta Z := Z - Z'$. Applying It\^o's lemma to $\| \delta Y \|^2$, we obtain
\[
	d \| \delta Y(t) \|^2 = 2 \langle \delta Y(t), d (\delta Y(t)) \rangle + \| \delta Z (t)\|^2 dt.
\]
We recall that the Hamiltonian is given by
\[
	H(t,x,u,y,z) := b'(t,x,u)y + \Tr(\sigma'(t,x,u)z) + L(t,x,u).
\]
We define an auxiliary function 
\[
\begin{split}
	f(s,x,u,z) &:= \nabla_x H(s,x,u,y,z) - \nabla_x b(s,x,u) y \\
			& = \nabla_x L(s,x,u) + \Tr \big(\nabla_x \sigma (s,x,u)z \big)
%	K^{Y,Z}_s := \nabla_x H(s,X^u_s,u_s,Y(s),Z(s)) - \nabla_x b(s,X^u_s,u_s) Y(s)
\end{split}
\]
to represent the ``undiscounted'' part of the driver term in (\ref{IHBSDE}). We emphasise that $f$ \textit{does not depend} on $y$. Then, by Girsanov's theorem, one can show that there exists a probability measure $\bQ \sim \bP$ such that the process
\[
\begin{split}
	M_t &:= \int_t^T \langle \delta Y(s),\delta Z(s) \rangle dW_s 
	\\ & \qquad + \int_t^T \langle \delta Y(s) , (f(s,X^u_s,u_s,Z(s)) - f(s,X^u_s,u_s,Z'(s))) \rangle ds
\end{split}
\]
is a $\bQ$-martingale. Thus, for all $0 \leq s \leq t \leq T$, we see that
\[
\begin{split}
	\bE_s^Q \| \delta Y(t) \|^2 & = \| \delta Y(s) \|^2  + \bE^Q_s \int_s^t \| \delta Z (l)\|^2 dl
	 \\ & \quad - 2\bE^Q_s \int_s^t \langle \nabla_x b(l,X^u_l,u_l) \delta Y(l),\delta Y(l) \rangle dl
\end{split}
\]
where $\bE^Q_s$ denotes the conditional expectation under $\bQ$ given $\F_s$. Using Assumption \ref{assumption_gradient}, we know that 
\[
	- 2 \langle \nabla_x b(l,X^u_l,u_l) \delta Y(l),\delta Y(l) \rangle  \geq 2 k \delta \| Y(l) \|^2,
\]
and hence we immediately see that 
\[
	\frac{d}{dt}\bE_s^Q \| \delta Y(t) \|^2 \geq 2k \bE_s^Q \| \delta Y(t) \|^2,
\]
leading to 
\[
	\bE_s^Q \| \delta Y(t) \|^2 \geq \| \delta Y(s) \|^2 e ^ {2k(t-s)},
\]
and hence 
\begin{equation}
	\| \delta Y(s) \|^2 \leq \bE_s^Q \| \delta Y(t) \|^2 e^{-2k(t-s)} \leq 2C' e^{-2k(t-s)}, 
\label{YYY}
\end{equation}
where $\bar{C}$ is the bound on $\| Y \|$. The right hand side in (\ref{YYY}) is independent of $T$ and collapses as $t\to\infty$. Hence $|\delta Y_s|=0$, from which we see $Y_s=Y'_s$ a.s. for every $s$, and hence $Y=Y'$ up to indistinguishability as $Y$ and $Y'$ are c\`adl\`ag.

We now show that a bounded solution exists. We first notice that there indeed exists a unique solution to the $T$-horizon BSDE (\ref{FHBSDE}) (see, for example, \cite{Royer}). In order to prove that $Y^T$ is bounded, we proceed similarly to above. Applying It\^o's lemma to $\| Y^T \|^2$, we obtain
\[
	d \|  Y^T(t) \|^2 = 2 \langle  Y^T(t), d  Y^T(t) \rangle + \| Z^T (t)\|^2 dt.
\]
As before, we know that there exists a probability measure $\bQ \sim \bP$ such that the process
\[
\begin{split}
	M_t &:= \int_t^T \langle Y^T(s), Z^T(s) \rangle dW_s \\ &\qquad+ \int_t^T \langle Y^T(s) , (f(s,X^u_s,u_s,Z^T(s)) - f(s,X^u_s,u_s,0)) \rangle ds
\end{split}
\]
is a martingale. Thus, taking conditional expectations under $\bQ$, for all $0 \leq s \leq t \leq T$, we see that
\[
\begin{split}
	\frac{d}{dt}\bE_s^Q \|  Y^T_t \|^2 &\geq -2 \bE \langle K^{Y^T,0}_t , Y^T_t \rangle + 2k \bE_s^Q \| \delta Y_t \|^2 \\
							& \geq -2 \bE C \| Y^T_t \| +  2k \bE_s^Q \| \delta Y_t \|^2
\end{split}
\]
as a simple application of Cauchy--Schwarz, where $C$ is the uniform bound on $\| f(s,X^u_s,u_s,0) \| =  \| \nabla_x L(t,X^u_t,u_t) \|$ from Assumption \ref{assumption_gradient}. It is clear that, for any $\epsilon > 0$, there exists a constant $C^{\epsilon} > 0$ such that $Cx \leq \epsilon x^2 + C^{\epsilon}$ for all $x \geq 0$. Using this fact, we arrive at 
\[
	\frac{d}{dt}\bigg[ \bE_s^Q \|  Y^T_t \|^2 - \frac{C^{\epsilon}}{k - \epsilon} \bigg] \geq 2(k-\epsilon) \bigg[ \bE_s^Q \|  Y^T_t \|^2 - \frac{C^{\epsilon}}{k - \epsilon} \bigg].
\]
By Gr\"onwall's lemma we see that 
\[
	 \bE_s^Q \|  Y^T_t \|^2 - \frac{C^{\epsilon}}{k - \epsilon}  \geq e^{2(k-\epsilon)(t-s)}\bigg[  \|  Y^T_s \|^2 - \frac{C^{\epsilon}}{k - \epsilon}   \bigg].
\]
Choosing $\epsilon < k$, and taking $t = T$, we conclude that 
\[
	\|  Y^T_s \|^2 \leq -e^{-2(k-\epsilon)(T-s)} \frac{C^{\epsilon}}{k - \epsilon} + \frac{C^{\epsilon}}{k - \epsilon} \leq \frac{C^{\epsilon}}{k - \epsilon},
\]
and thus $Y^T$ is uniformly bounded. Denoting $\delta Y := Y^T - Y^{T'}$, $\delta Z := Z^T - Z^{T'}$, we proceed exactly as we did to prove uniqueness, to arrive at
\begin{equation}
	\| \delta Y_s \|^2 \leq \bE_s^Q \| \delta Y_t \|^2 e^{-2k(t-s)} \leq \frac{2 C^{\epsilon}}{k-\epsilon} e^{-2k(t-s)}, 
\label{YYYY}
\end{equation}
where $\bQ \sim \bP$ is the probability measure such that the process 
\[
\begin{split}
	M_t &:= \int_t^T \langle \delta Y(s),\delta Z_s \rangle dW_s 
		\\ & \qquad + \int_t^T \langle \delta Y(s) , (f(s,X^u_s,u_s,Z^T(s)) - f(s,X^u_s,u_s,Z^{T'}(s))) \rangle ds
\end{split}
\]
is a martingale for $0 \leq t \leq T \leq T'$. From (\ref{YYYY}) we see that $Y^T_t$ is a Cauchy sequence in $T$. Therefore the limit exists, and we denote it $Y_t$. The uniform bound also holds for $Y_t$, and convergence uniformly on compacts is clear from ($\dagger$). In order to see that $Z^T$ is Cauchy as well, we apply It\^o's formula to $\| \delta Y \|^2$, and take expectations under $\bQ$. Then 
\[
	\bE^Q \int_0^t \| \delta Z (s)\|^2 ds \leq -2k \bE^Q \int_0^t \| \delta Y(s) \|^2 + \bE \| \delta Y(t) \|^2 - \| \delta Y(0) \|^2,
\]
and, given (\ref{YYYY}), the claim follows. Therefore, the limit as $T \to \infty$ exists for the sequence $\{ Z^T_t \}$. Taking $Z$ as its limit we obtain our desired solution $(Y,Z)$. 

\qed

\noindent Having proven the existence of a unique solution to the adjoint BSDE on an infinite horizon, we aim to establish a relation to ergodicity of $\{X^u_t\}_{t \geq 0}$ under all controls $u$. We begin with the following lemma:

\begin{lemma} Suppose for all $(t,u) \in \bR^+ \times \U$, the function $b(t,\cdot,u) \in C^1$. Then the following conditions are equivalent: 
\begin{enumerate}[(i)]
\item The function $b(t,\cdot,u)$ is strongly dissipative, that is 
\[
	\langle b(t,x,u) - b(t,y,u), x-y \rangle \leq -k \| x-y \|^2
\]
holds for all $t \geq 0$, $x,y \in \bR^n$ and for all $u \in \U$. 
\item The gradient matrix $\nabla_x b(t,\cdot,u)$ is negative definite, that is 
\[
	\langle \nabla_x b(t,x,u)y,y \rangle \leq -k \| y \|^2 
\]
holds for all $t \geq 0$, $x,y \in \bR^n$ and for all $u \in \U$. 
\end{enumerate}
\label{L1}
\end{lemma}
\noindent \textbf{Proof:} \textbf{(i) $\to$ (ii):} Given any $\alpha \in [0,1]$, let $x(t) := \alpha x + (1-\alpha)y$. Define the function $g:[0,1] \to \bR$ by
\[
	g(\alpha) = \langle b(t,x(\alpha),u) - b(t,y,u), x-y \rangle - k (1-\alpha) \| x-y \|^2 ,
\] 
where $k$ is the same as in Assumption \ref{assumption_gradient}. We notice that $g(0) = -k \| x-y \|^2$, and $g(1) = \langle b(t,x,u) - b(t,y,u), x-y \rangle$. We see also that 
\[
	g'(\alpha) = \langle \nabla_x b(t,x(\alpha),u)(x-y),x-y \rangle + k \| x-y \|^2 ,
\]
and therefore $g'(\alpha) \leq 0$ for all $t \in [0,1]$. Therefore $g(1) \leq g(0)$, finishing the proof.

\noindent \textbf{(ii) $\to$ (i):} By assumption, we know that, for all $y \in \bR^n$ and for all $\epsilon > 0$,
\[
	\langle b(x+\epsilon y) - b(y), \epsilon y \rangle \leq -k \epsilon^2 \| y \|^2.
\]
Dividing both sides by $\epsilon^2$, and taking the limit, we arrive at 
\[
		\bigg\langle \lim_{\epsilon \to 0} \frac{b(x+\epsilon y) - b(y)}{\epsilon},  y \bigg\rangle \leq -k \| y \|^2,
\]
and the result follows, given that the function $b$ is continuously differentiable.

\qed

\noindent We have thus proven that the conditions proposed for the adjoint BSDE to admit a solution are equivalent to strong dissipativity of the forward process. In the next section we will show that this in turn implies ergodicity, in the sense that the laws of two solutions to the equation (\ref{CDX}) with different initial values get exponentially close in time. This will allow us to establish the connection to the representation via ergodic BSDEs. 

\begin{remark} One could wish to establish connection between the solvability of the adjoint BSDE in infinite horizon and weak dissipativity of the forward process. The reason being that it is the largest class for which we can expect to generally prove ergodicity. However, even in one dimension, the condition 
\[
	\langle \nabla_x b(t,x,u)y,y \rangle \leq -k \| y \|^2 	
\]
cannot be relaxed without losing the existence of solutions to the adjoint equation on infinite horizons. By Lemma \ref{L1} we then know that $b$ is strongly dissipative. 
\end{remark}

\section{Structural properties of the forward process}

\noindent The main goal of this section is to show that, if Assumptions \ref{existence_solution} and \ref{assumption_gradient} are satisfied, then, for each $\{ u_t \}_{t \geq 0}$, the process $X^u$ is strong Feller and irreducible. We then proceed to prove that the laws of two processes satisfying (\ref{CDX}) and started at different points get exponentially close as $t \to \infty$. 

\begin{remark} In the sequel we omit $u$ from the dynamics of $X$ for notational simplicity.
\end{remark}

\noindent We start with an auxiliary lemma that will be crucial in showing the exponential convergence of laws. 

\begin{lemma} Let $\{X^x_t\}_{t \geq 0}$ be a solution to (\ref{CDX}) with $X^x_0 = x$. Suppose also that Assumption \ref{assumption_gradient} is satisfied. Then there exist constants $c \geq 0$ and $\mu > 0$ such that 
\[
	\bE \| X^x_t \|^2 \leq  \| x \|^2 e^{-\mu t} + c
\]
for all $t \geq 0$.
\label{L5}
\end{lemma}

\noindent \textbf{Proof:} Applying It\^o's formula to $\| X^x_t\|^2$ and taking expectations, we obtain
\[
	\bE \| X^x_t \|^2 = \|x\|^2 + \bE \int_0^t G(s)ds,
\]
where
\[
\begin{split}
	G(s) &= 2\langle b(s,X^x_s), X^x_s  \rangle + \| \sigma(s,X^x_s) \|^2 \\
		& \leq -2k \| X^x_s\|^2 + \langle b(s,0) X^x_s \rangle + \| \sigma(s,X^x_s) \|^2.
\end{split}
\]
Taking into account Assumption \ref{drift} and applying Cauchy--Schwarz, we see that for any $\epsilon > 0$, there exists $C^{\epsilon}$, such that 
\[
	G(s) \leq -2(k - \epsilon)\| X^x_s \|^2 + C^{\epsilon},
\]
for $0 \leq s \leq t$. Then 
\[
	\frac{d}{dt} \bE \| X^x_t \|^2 \leq -(k - \epsilon)\bE \| X^x_t \|^2 + C^{\epsilon}, \quad t \geq 0.
\]
We can pick any $\epsilon < k$, and apply Gr\"onwall's lemma to arrive at 
\[
	\bE \| X^x_t \|^2 \leq \| x \|^2 e^ {-(k-\epsilon)t} + \frac{C^{\epsilon}}{k-\epsilon},
\]
concluding the proof.

\qed

\subsection{Strong Feller Property} 
We begin by defining a number of objects we will be repeatedly using in the sequel.  
\begin{definition} Here and for the rest of the paper, we use the standard notation 
\[
	P_t \psi (x) := \bE \big[ \psi (X^x_t) \big]
\]
for the transition operator (or the semigroup associated with the process $\{X^x_t\}_{t \geq 0}$), and $\| \psi \|_0 = \sup_{u \in \bR^n}\| \psi(u) \|$. We also write 
\[
	P_t(x,A) = P_t [\bone_{A}] (x), \quad A \in \bR^n
\]
for the law of the process $\{X^x_t\}_{t \geq 0}$.
\end{definition}
\noindent We first state the main result of this section, the proof of which will require two auxiliary lemmas:
\begin{theorem} Suppose Assumptions \ref{existence_solution} and \ref{assumption_gradient} hold. Then for all $\psi \in \B_b(\bR^n)$ and $t > 0$,
\begin{equation}
	| P_t \psi (x) - P_t \psi (y) | \leq C_t \| x - y \| \| \psi(u) \|_0
\label{Feller_estimate}
\end{equation}
holds for some constant $C_t$.
\label{Feller}
\end{theorem}

\begin{lemma} Let $c>0$ and $t >0$ be fixed. Then the following conditions are equivalent:
\begin{enumerate}
\item For all $\psi \in C^2_b(\bR^n)$ and for all $x,y \in \bR^n$ we have
\[
	| P_t \psi (x) - P_t \psi (y) | \leq c \| x - y \| \| \psi(u) \|_0.	
\]
\item For all $\psi \in B_b(\bR^n)$ and for all $x,y \in \bR^n$ we have
\[
	| P_t \psi (x) - P_t \psi (y) | \leq c \| x - y \| \| \psi(u) \|_0.	
\]
\item For all $x,y \in \bR^n$
\[
	| P_t(x,\cdot) - P_t(y,\cdot) |_{TV} \leq c \| x-y \|.
\]
\end{enumerate}
\label{L4}
\end{lemma}
\noindent \textbf{Proof:} See Lemma 9.36 in \cite{PRATO}. 

\qed

\noindent In addition to the stochastic differential equation (\ref{CDX}), we now introduce a corresponding velocity process (a directional derivative with respect to the initial value) $t \to V^h(t), t \geq 0$, defined by 
\[
	V^h(t) := \langle D_x X^x_t , h\rangle
\]
for $h \in \bR^n$ and $t \geq 0$. We notice that this is nothing but a directional derivative of the stochastic flow $x \mapsto X^x$ defined by the SDE (\ref{CDX}). We note that under our regularity assumptions on the coefficients $b$ and $\sigma$ this is a meaningful object (see Theorem 9.6 in \cite{PRATO} for details). We now prove a crucial estimate on the norm of $V^h(t)$.
\begin{lemma} Suppose Assumption \ref{assumption_gradient} holds. Then, for all $t \geq 0$, there exists a constant $c_t$, such that
\[
	\bE \| V^h(t) \|^2 \leq c_t \| h \|^2,
\]
for all $h \in \bR^n$.
\label{L3}
\end{lemma}
\noindent \textbf{Proof:} By definition, the process $X^x$ satisfies the integral version of (\ref{CDX}), that is 
\[
	X^x_t = x + \int_0^t b(s,X^x_s)ds + \int_0^t \sigma(s,X^x_s)dW_s.
\]
taking directional derivatives of both sides, we obtain 
\[
	\langle D_x X^x_t, h \rangle =  h + \int_0^t \nabla_x b(s,X^x_s) \langle D_x X^x_s, h \rangle ds + \int_0^t \nabla_x \sigma(s,X^x_s)  \langle D_x X^x_s, h \rangle dW_s,
\]
which can be immediately rewritten as 
\[
	V^h(t)=  h + \int_0^t \nabla_x b(s,X^x_s) V^h(s) ds + \int_0^t \nabla_x \sigma(s,X^x_s)  V^h(s) dW_s.
\]
Due to the Assumption \ref{existence_solution}, we see that $\| \nabla_x \sigma(s, X^x_s) h \|^2 \leq \omega \| h \|^2$, where $\omega$ is the Lipschitz constant for $\sigma(t,\cdot)$. Applying It\^o's lemma to $\| V^h(t) \|^2$, taking expectations and using Assumption \ref{assumption_gradient}, we arrive at 
\[
	\bE\| V^h(t) \|^2 \leq \| h \|^2 + (\omega-k) \int_0^t \bE\| V^h(s) \|^2 ds, 
\]
and the conclusion follows by Gr\"onwall's lemma with $c_t = e^{(\omega-k)t}$. 

\qed

\noindent In order to proceed, we need a version of the Bismut--Elworthy formula:

\begin{lemma} In the notation above, the following holds:
\begin{equation}
	\bE \bigg[  \psi(X^x_t) \int_0^t  \sigma^{-1}(s,X^x_s)V^h(s) dW_s \bigg] = t \langle h, D_x P_t \psi (x) \rangle, 
\label{XXX}
\end{equation}
where $\psi \in C^2_b(\bR^n)$. 
\end{lemma}
\noindent \textbf{Proof:} See, for example, \cite{JVC}. 

\qed

\noindent \textbf{Proof: (Theorem \ref{Feller}).} Using formula (\ref{XXX}), we deduce 
\[
	\| \langle h, D_x P_t \psi (x) \rangle \| ^2 \leq \frac{1}{t^2} \| \psi(u) \|_0^2 \bE \bigg[ \int_0^t \| \sigma^{-1}(s,X^x_s)V^h(s) \|^2 ds \bigg],
\]
and, using Lemma \ref{L3}, we obtain the following estimate:
\[
	\| \langle h, D_x P_t \psi (x) \rangle \| ^2 \leq \frac{\underline{\sigma}^{-2}}{t^2} \| \psi(u) \|_0^2 \| h \|^2 \int_0^t c_s ds,
\]
and therefore
\[
	\| D_xP_t \psi(x) \| \leq C_t := \frac{\underline{\sigma}^{-1}}{t} \| \psi(u) \|_0  \bigg(\int_0^t c_s ds\bigg)^{1/2},
\]
and thus 
\[
	| P_t \psi (x) - P_t \psi (y) | \leq C_t \| x - y \| \| \psi(u) \|_0
\]
holds for any $x,y \in \bR^n$, $t \geq 0$ and $\psi \in C^2_b(\bR^n)$. The conclusion follows immediately by Lemma \ref{L4}.

\qed

\begin{remark} It is straightforward to see that all the estimates established above would still hold if we replaced $P_t \psi (x)$ with $P(\tau,t)[f](x) = \bE \psi(X^{x,\tau}_t)$, where  $X^{x,\tau}$ is the solution to 
\begin{equation}
	dX^{x,\tau}_t = b(t,X^{x,\tau}_t)dt + \sigma(t,X^{x,\tau}_t)dW_t, \quad X_{\tau} = x.
\label{ISDE}
\end{equation}
In this case all time dependent constants $C_t$ will become $C_{t-\tau}$. 
\label{R1}
\end{remark}

\subsection{Irreducibility}
The goal of this section is to prove the irreducibility of the forward process, i.e. that it at any time it can be in any open ball with positive probability. The main result is as follows:
\begin{theorem} Let $\{X^x_t\}_{t \geq 0}$ be the solution to (\ref{CDX}), and let $B_r(x)$ denote an open ball of radius $r$ around $x \in \bR^n$. Then 
\[
	\bP(X^x_t \in B_{r}(z)) > 0
\]
for any $t > 0$, $z \in \bR^n$, $r > 0$.
\label{irreducibility}
\end{theorem}
\noindent \textbf{Proof:} We first notice that irreducibility for $X^x$ is equivalent to irreducibility of the (strong) solution to 
\[
	dY_t = (b(t,Y_t) + \phi(t))dt + \sigma(t,Y_t)dW_t, \quad Y_0 = x,	
\]
where $\phi(\cdot,\cdot)$ is a bounded function. In order to see this, we recall that $\sigma^{-1}(t,x)$ is bounded away from zero, and therefore, for any $T > 0$, we can find a probability measure $\bQ^T \sim \bP$, such that under $\bQ^T$, the process $dW^{Q}_t = dW_t + \phi(t)dt$ defines a Brownian motion for all $t \in [0,T]$. Since these measures are equivalent, so are the laws of $X^x$ and $Y$ on $[0,T]$. For details see Lemma 7.3.2 in \cite{G.DaPrato1996}. The rest of the proof is inspired by Lemma 7.3.3 in \cite{G.DaPrato1996}. We fix $t > 0$, $x,z \in \bR^n$, $r > 0$. We will now construct a bounded function $\phi$, such that the corresponding solution satisfies 
\[
	\bP(\| Y_t - z \| < r) > 0.
\]
By the reasoning above this is enough to prove the claim. Let $t_1 \in [0,t)$ be a moment of time to be chosen later. We then define $\phi$ in the following way:
\[
	\phi(s) := \begin{cases} 0, \quad &s \in [0,t_1), \\
					    \xi(s,Y_{t_1}), \quad &s \in [t_1,t],	
	\end{cases}
\]	
where
\[
	\xi(s,u) = \begin{cases} \frac{y - b(s,0) - u}{t-t_1}, \quad  &\| u \| < R, \\
						0, \quad &\text{otherwise},
	\end{cases}
\]
where $y,R > 0$ are constants to be chosen later. We start by picking any $y$ such that $\| y - a \| < r/3$. Then, if $\| Y_{t_1} \| < R$, we have 
\[
	Y_t = y + \int_{t_1}^t [b(s,Y_s) - b(s,0)]ds + \int_{t_1}^t \sigma(s,Y_s)dW_s =: y + I,
\]
and thus 
\[
	\bP\bigg(\| Y_t - z \| < r\bigg) \geq \bP \bigg(\| Y_{t_1} \| < R, \| I \| \leq 2r/3\bigg) \geq \bP \bigg(\| Y_{t_1} \| < R\bigg) - \bP\bigg(\| I \| > 2r/3\bigg),
\]
and it remains to pick $t_1$ such that 
\[
	\bP(\| I \| > 2r/3) \leq 1/4,
\]
and $R >0$ such that 
\[
	\bP (\| Y_{t_1} \| < R) \geq 3/4.
\]
We notice that the estimate in Lemma \ref{L5} is true for the process $Y$ as well. Therefore, taking into account that the function $b(s,\cdot)$ is Lipschitz, and applying Markov's inequality, we have
\[
	\bP(\| I \| > 2r/3) \leq \frac{3}{2r}\bE \big[ \| I \| \big] \leq \frac{3\bigg( (t-t_1) K \big[\|x\| + c \big]  + \sqrt{t-t_1}\bar{\sigma}\bigg)}{2r},
\]
and hence the choice of $t_1$ is always possible. The choice of $R$ is a direct application of Lemma \ref{L5}. 

\qed

\subsection{Main result}

We are now ready to prove the exponential convergence of laws for two solutions of the SDE (\ref{CDX}) with different initial values. This result makes it possible to talk about Ergodic BSDEs in the present framework, and gives us hope to relate two existing approaches to Ergodic Optimal Control. 
\begin{theorem} Suppose Assumptions \ref{existence_solution} and \ref{assumption_gradient} hold. Then there exist constants $C > 0$ and $\rho > 0$ such that, for any bounded continuous function $\psi : \bR^n \to \bR$,
\begin{equation}
	|P(\tau,t)[\psi](x)  - P(\tau,t)[\psi](y) | \leq C (1 + ||x||^2 + ||y||^2)e^{-\rho (t-\tau)}\|\psi(u)\|_0,
\label{Estimate_main}	
\end{equation}
where $P(\tau,t)[f](x)$ is defined as in Remark \ref{R1}.
\label{mainThm}
\end{theorem}
\noindent \textbf{Proof:} We will not provide a detailed proof, since most of it is would closely follow, e.g. \cite{Levy_paper}, but we will outline a strategy. The main idea is to construct a discrete time Markov chain $(\bar{X}_{k\tT}, \bar{Y}_{k\tT})$ with a time step $\tT$ on the space $\bR^{2n}$, whose transition operator is determined by the law of the processes $\{X^x_t\}_{t \geq 0}$ and $\{X^y_t\}_{t \geq 0}$. We then prove the necessary convergence of laws for this chain. 

\begin{itemize}
\item \textbf{(Step 1): } The first step is almost identical to the corresponding one in the proof of Theorem 5 in \cite{Levy_paper}, except for a new proof of Lemma \ref{L5}. We start by showing that we can choose a time step $\tT > 0$ and a radius $R > 0$, such that, if we observe two independent solution processes $X^x$ and $X^y$ only at times $\{ n\tT \}_{n \in \bN}$, there is an \textit{exponential bound} on the waiting time for both $X^y_{n\tT}$ and $X^x_{n\tT}$ to enter $B_R(0)$. The independence here is understood in the the sense that we take two independent copies ($\tilde{W}$ and $\bar{W}$) of the Brownian motion $W$. We set 
\[
	\bar{X}_{k\tT} = X^x_{k\tT}, \quad \bar{Y}_{k\tT} = X^y_{k\tT},
\]
for $ k \leq \tau := \inf\{ n: X^x_{n\tT} \in B_R(0) , X^y_{n\tT} \in B_R(0)\}$. 
\item \textbf{(Step 2):} We will explain this part in a little more detail, since the method here is different from that in \cite{Levy_paper}. Once $X^x_{k\tT}$ and $X^y_{k\tT}$ are in $B_R(0)$ for some $k \geq 0$, we lift the independence assumption and construct two solutions $X^{X^x_{k\tT}, k\tT}$ and $X^{X^x_{k\tT}, k\tT}$ to (\ref{ISDE}) on $[k\tT, (k+1)\tT]$ with initial conditions $X^x_{k\tT}$ and $X^y_{k\tT}$ respectively. By irreducibility of solutions to (\ref{ISDE}) and compactness of $B_R(0)$, for any fixed $\epsilon > 0$, there exists a constant $\delta_{\epsilon} > 0$ such that 
\[
	\bP \bigg(  X^{X^x_{k\tT}, k\tT}_{(k+0.5)\tT} \in B_{\epsilon}(0),  X^{X^y_{k\tT}, k\tT}_{(k+0.5)\tT} \in B_{\epsilon}(0)  \bigg) > \delta_{\epsilon}.
\]
We proceed as follows: we know (see Theorem 5.2 in \cite{Lindv}) that there exists a coupling $(Z,Z')$ such that 
\[
	\bP(Z \neq Z') = \big\| P_{(k+1)\tT}(X^x_{(k+0.5)\tT},\cdot)  - P_{(k+1)\tT}(X^y_{(k+0.5)\tT},\cdot)\big\|_{TV},
\]
where $P_{(k+1)\tT}(x,\cdot)$ is the law of the solution to (\ref{ISDE}) with $X_{(k+0.5)\tT} = x$, and $Z,Z'$ are independent conditioned on the event $\{Z \neq Z'\}$. By Theorem \ref{Feller} and Lemma \ref{L4}, there exists $\epsilon > 0$ such that, for $x,y \in B_{\epsilon}(0)$,
\[
	\big\| P_{(k+1)\tT}(X^x_{(k+0.5)\tT},\cdot)  - P_{(k+1)\tT}(X^y_{(k+0.5)\tT},\cdot)\big\|_{TV} \leq \frac{1}{2},
\]
and therefore 
\[
	\bP(Z \neq Z') \geq \frac{1}{2}.
\]
Now we set 
\[
	\bar{X}_{(k+1)\tT} = \begin{cases} Z, \quad &(X^x_{(k+0.5)\tT},X^y_{(k+0.5)\tT}) \in B_{\epsilon}(0) \times B_{\epsilon}(0),
	 \\ X^x_{(k+1)\tT}, \quad &\text{otherwise}. \end{cases}
\]
Define $\bar{Y}_{(k+1)\tT}$ similarly. We then infer that, for the constructed solutions, 
\[
	\bP(\bar{X}_{(k+1)\tT} = \bar{Y}_{(k+1)\tT}) \geq \frac{\delta_{\epsilon}}{2},
\]
 which is a bound from below uniformly in $k$. 
\item \textbf{(Step 3): } We then iterate these arguments to show that the probability that the two processes we are constructing have not met decays exponentially in time.  

\end{itemize}

\qed

%\begin{remark} We have shown that under the conditions required for the existence of a solution to the adjoint BSDE (\ref{adjoint}), the forward process is ergodic under all controls. This is v
%\end{remark}

\section{Connection to Ergodic BSDEs}

In this section we explore an alternative way of looking at optimal ergodic control problems. Recall that we aim to minimise 
\begin{equation}
	J(x_0, u) = {\lim\sup}_{T\to\infty} T^{-1} \bE \bigg[\int_0^T L(t,X^u_t, u_t) dt\bigg],
\label{J}
\end{equation}
over the space $\U$ of controls, where $L$ is uniformly bounded by some constant $C > 0$. If we assume that we can attack this problem though the adjoint BSDE approach, by Lemma \ref{L1} the process $X^u_t$ is ergodic under all control processes $\{u_t\}$. In order to ensure the existence of an ergodic measure, we impose the following conditions:
\begin{assumption}There exists a positive period $T^*$, such that 
\begin{enumerate}[(i)]
\item For all admissible control processes $\{u_t\}_{t \geq 0}$, $u_{t + T^*} = u_t$.
\item The coefficients $b$ and $\sigma$ in the dynamics of $X$ are $T^*$-periodic.
\item The cost function $L$ is $T^*$-periodic.
\end{enumerate}
\end{assumption}
\noindent We know (see, e.g. \cite{Knable}), that there exists a unique invariant measure $\mu$ for the semigroup associated with the process $X^u$, such that there exists ergodic average $\lambda^u$, satisfying 
\[
	\bar{\lambda^u} = \limsup_{T \to \infty}\bE \int_0^T  L(t,X^{u}_t,u_t)dt = \int_{[0,T^*] \times \bR^n} L(t,x)\mu(dt,dx).
\]
\noindent As a special case of a much more general result (see, e.g. Theorem 10 in \cite{Levy_paper}) one could also prove the following:

\begin{theorem} For any fixed control process $\{u_t\}_{t \geq 0}$, there exists a solution triple $(Y^u,Z^u,\lambda^u)$ to the equation 
\begin{equation}
	 Y^u_t = Y^u_T+\int_t^T[L(s,X^u_t, u_s)-\lambda^u] du - \int_t^T Z^u_s dW_s.
\label{EBSDE}
\end{equation}
There also exists a deterministic function $v^u(t,x)$, such that $Y^u_t = v^u(t,X^u_t)$, $v^u$ is locally Lipschitz in space and globally in time, and 
\[
	\| v^u(t,x) \| \leq C' (1 + \| x \|^2)
\]
for some constant $C' > 0$. 
\end{theorem}

\begin{remark} Note that in the present context, the driver function $L$ is time dependent. This fact poses no problem to the construction of a solution via vanishing discount approach, since the existence of a solution to the discounted infinite horizon BSDE can be established in the most general case (see, e.g. Theorem 8 in \cite{Levy_paper}). One could also still show, that if there is another Markovian solution $(\tilde{v}^u,\tilde{Z}^u,\tilde{\lambda}^u)$, such that $\tilde{v}^u$ is of polynomial growth in $x$, then $\lambda^u = \tilde{\lambda}^u$. \end{remark}

\noindent We now demonstrate that for a fixed control process $u$, $\lambda^u$ is the corresponding ergodic average. In other words,
\[
	\lambda^u = \bar{\lambda}^u.
\]
In order to see this, denote $(Y^{u,T}, Z^{u,T})$ to be a unique solution to the $T$-horizon BSDE 
\[
	 Y^{u,T}_t = \int_t^T L(s,X^u_t, u_s) du - \int_t^T Z^{u,T}_s dW_s.	
\]
It is clear that 
\[
	Y^{u,T}_0 = \bE \int_0^T L(t,X^u_t, u_t) dt.
\]
Comparing the dynamics of $Y^u$ and $Y^{u,T}$ and taking expectations, we immediately see that  
\[
	Y^{u,T}_0 - v^u(0,x_0) = - \bE v^u(T,X^u_T) + \lambda T.
\]
By polynomial growth we know that 
\[
	\| v^u(0,x_0) \| \leq C' (1 + \| x_0 \|^2), \quad \| v^u(T,X^u_T) \| \leq C' (1 + \| X^u_T \|^2).
\]
Recall (see Theorem 4, \cite{Levy_paper}) that 
\[
	\bE \| X^u_T \|^2 \leq C''(1 + \|x_0 \|^2)
\]
for some constant $C'' > 0$ independent of time. Therefore 
\[
	J(x_0,u) = \lim\sup_{T \to \infty} \frac{Y^{u,T}}{T} = \lambda^u.
\]
Then the optimal control process $\{u^*_t\}$ satisfies 
\[
	\lambda^{u^*} = \inf_{u} \lambda^u.
\]
We note that since we control both drift and volatility of the forward process, it can not be represented as measure change, and therefore we are not able to arrive at the optimal value by solving just one ergodic BSDE with infimum of some Hamiltonian in the driver.

\section{Conclusion}

We have established a duality between the conditions implying the ergodicity of the forward process and the conditions required for the stochastic maximum principle approach to the problem of optimisation to make sense. The result confirms our intuition: one would expect that in order to have any hope for solving an ergodic control problem via probabilistic techniques, the ergodicity of the controlled process needs to be assumed. The converse is also true, in a sense that if the forward process is strongly dissipative under all controls, then we can deal with the problem via the adjoint equation technique. We have shown that this is indeed the case. One direction for future work would be to find a different approach that would allow us to attack the case of weakly dissipative drifts within the strong formulation.

\end{document}